\documentclass[aop,seceqn,MSNbibl,citesort,dvips]{arximspdf}

\doi{10.1214/10-AOP638}
\volume{40}
\issue{2}
\pubyear{2012}
\firstpage{765}
\lastpage{787}

\makeatletter

\renewcommand{\overline}{\bar}

\newtheorem{theorem}{Theorem}[section]
\newtheorem{prop}[theorem]{Proposition}
\newtheorem{lemma}[theorem]{Lemma}
\newtheorem{cor}[theorem]{Corollary}

\newproclaim{remark}[theorem]{Remark}

\newcommand{\R}{\mathbb{R}}
\newcommand{\Nn}{\mathbb{N}}

\newcommand{\C}{\mathbb{C}}
\newcommand{\Q}{\mathbb{Q}}
\newcommand{\E}{\mathbb{E}}
\newcommand{\El}{\mathrm{E}}
\newcommand{\Pl}{\mathrm{P}}
\newcommand{\tp}{{\overline{\psi}}}

\newcommand{\ta}{\tilde{\alpha}}
\newcommand{\F}{\mathcal{F}}
\newcommand{\I}{\mathcal{I}}
\newcommand{\Ip}{\mathcal{I}_{\psi}}
\newcommand{\Ipg}{\mathcal{I}_{\psi_{\gamma}}}
\newcommand{\It}{\I_{\psi_{\theta}}}

\newcommand{\ann}{a_n(\psi;\alpha)}

\newcommand{\km}{\kappa(a)}
\newcommand{\N}{\mathcal{N}}

\makeatother

\begin{document}
\begin{frontmatter}

\title{Law of the absorption time of some positive self-similar Markov processes\thanksref{T1}}
\runtitle{Absorption time of some Markov processes}

\thankstext{T1}{Supported in part by Swiss National Fund Grant
2000021--121901.}

\begin{aug}
\author[A]{\fnms{P.} \snm{Patie}\corref{}\ead[label=e1]{ppatie@ulb.ac.be}}
\runauthor{P. Patie}
\affiliation{Universit\'{e} Libre de Bruxelles}
\address[A]{D\'{e}partement de Math\'{e}matiques\\
Universit\'{e} Libre de Bruxelles\\
Boulevard du Triomphe\\
B-1050, Bruxelles\\
Belgique\\
\printead{e1}} 
\end{aug}

\received{\smonth{4} \syear{2010}}
\revised{\smonth{12} \syear{2010}}

%
\begin{abstract}
Let $X$ be a spectrally negative self-similar Markov process with $0$
as an absorbing state. In this paper, we show that the distribution of
the absorption time is absolutely continuous with an infinitely
continuously differentiable density. We provide a power series and a
contour integral representation of this density. Then, by means of
probabilistic arguments, we deduce some interesting analytical
properties satisfied by these functions, which include, for instance,
several types of hypergeometric functions. We also give several
characterizations of the Kesten's constant appearing in the study of
the asymptotic tail distribution of the absorbtion time. We end the
paper by detailing some known and new examples. In particular, we
offer an alternative proof of the recent result obtained by Bernyk,
Dalang and Peskir [\textit{Ann. Probab.} \textbf{36} (2008) 1777--1789]
regarding the law of the maximum of spectrally positive L\'{e}vy stable
processes.
\end{abstract}

%

\pdfkeywords{60E07, 60G18, 60G51, 33E30, Self-similar processes,
absorption time, Levy processes, exponential functional,
generalized hypergeometric functions}

\end{frontmatter}

\section{Introduction}

Let $X=((X_t)_{t\geq0}, (\Q_x)_{x>0})$ be a self-similar Hunt process
with values in $ [0,\infty)$. It means that $X$ is a right-continuous
strong Markov process with quasi-left continuous trajectories and there
exists $\alpha>0$ such that $X$ enjoys the following self-similarity
property: for each $c>0$ and $x\geq0$,
\[
\mbox{the law of the process } (c^{-1}X_{c^{\alpha}t})_{t\geq0},
\mbox{under } \Q_x, \mbox{ is } \Q_{x/c}.
\]
$1/\alpha$ is called the index of self-similarity. The purpose of the
paper is to describe the law of the stopping time
\[
T_0=\inf\{s>0;   X_s=0\}
\]
with the usual convention that $\inf\{\varnothing\}=\infty$. The class of
positive self-similar Markov processes (for short pssMp) has been
introduced and studied by Lamperti~\cite{Lamperti-72}. In particular,
he showed that for each fixed $\alpha> 0$, there is a~bijective correspondence
between pssMp with index $\alpha$ and (possibly killed) real-valued L\'{e}vy
processes, that is, processes with stationary and independent increments.
More specifically, by introducing the additive functional
\[
\Sigma_t=\inf\biggl\{s>0;  A_s =\int_0^sX_r^{-\alpha}\,dr>t\biggr\},
\]
Lamperti \cite{Lamperti-72} showed that the process $\xi=(\xi
_t)_{t\geq
0}$, defined by
%
\begin{equation} \label{eq:ss}
\xi_t = \log(X_{\Sigma_t}),\qquad 0\leq t<T_0,
\end{equation}
is a (possibly killed) L\'{e}vy process. We denote the law of the
process $\xi$ when starting at $0$ by $\mathrm{P}$.
It is plain that
\[
\Sigma_t = \int_0^t e^{\alpha\xi_s} \,ds
\]
and writing $q\geq0$ for the killing rate of the L\'{e}vy process, one
gets the identity in distribution
\[
(T_0,
\Q_x)\stackrel{(d)}{=}(x^{\alpha}\Sigma_{{\mathbf{e}}_{q}},\mathrm{P}),
\]
where ${\mathbf{e}}_q$ is an independent exponential random variable of parameter
$q$ (we have ${\mathbf{e}}_0=\infty$). Lamperti \cite{Lamperti-72} explained
that, either $q>0$ and $X$ reaches $0$ by a jump, that is,
\[
\Q_x(X_{T_{0-}}>0, T_0 <\infty)=1\qquad \forall x>0,
\]
or $\xi$ drifts to $-\infty$ and $X$ reaches $0$, that is,
\[
\Q_x(X_{T_{0-}}=0, T_0 <\infty)=1\qquad \forall x>0.
\]
We gather these two possibilities in the following hypothesis.
\begin{longlist}[H:]
\item[H:]\hypertarget{HypoH}\mbox{}
\textit{Either} $q>0$ \textit{or}
$\lim_{t\rightarrow\infty} \xi_t = -\infty$ \textit{a.s.}
\textit{and} $q=0$.
\end{longlist}

The law of $T_0$ or equivalently of $\Sigma_{{\mathbf{e}}_{q}}$ turns out
to be a
key object in various settings. It appears, for instance, in the study
of coagulation-fragmentation processes \cite{Bertoin-02-f} and
continuous state branching processes with immigration \cite
{Patie-CBI-09}. We also mention that, recently, in the SLE context,
Alberts and Sheffield~\cite{Alberts-Sheffield-08} describe a
measure-valued function supported on the intersection of a~chordal
SLE$(\kappa)$ curve with $\R$, $4<\kappa<8$, in terms of the law of the
absorption time~$T_0$ of some Bessel processes which form the class of
pssMp having continuous trajectories.
The law of $\Sigma_{{\mathbf{e}}_{q}}$ is also critical for the pricing
of Asian
options in mathematical finance (see, e.g., \cite
{Patie-Asian-09}), but also for computing perpetuities in insurance
mathematics (see, e.g., \cite{Dufresne-90}).

Unfortunately, beside some isolated cases the distribution of $T_0$ is
not attainable. We mention the papers
\cite{Carmona-Petit-Yor-97,Gjessing-Paulsen-97} and \cite{Patie-CBI-09} where such
examples can be found and refer to the survey paper \cite{Bertoin-Yor-05} for a
description of these cases. Besides, two notable exceptions might be
worth mentioning: when $X$ is a Bessel process of negative\vadjust{\goodbreak} index and
when $X$ is a regular spectrally negative stable L\'{e}vy process
killed upon entering the negative half-line. In the former case,
several proofs can be found in the literature, see, for instance, the
excellent monograph of Yor \cite{Yor-01} and the more recent survey
papers of Matsumoto and Yor \cite{Matsumoto-Yor-05-1} and \cite
{Matsumoto-Yor-05-2}. However, most of the proofs rely on the knowledge
of the semigroup of Bessel processes. For the second case,
Bernyk, Dalang and Peskir \cite{Bernyk-Dalang-Peskir-08} derive a representation of
the distribution of $T_0$ by inverting, in a nontrivial way, the known
expression of the Wiener--Hopf factorization of stable one-sided L\'
{e}vy processes. Our approach will differ from these two cases since we
do not have, in general, access neither to the semigroup of $X$ nor to
the Laplace transform of $T_0$.

The remaining part of the paper is organized as follows. In the next
section, we state our main results including the smoothness and the
representation as an absolutely convergent power series of the
distribution of $T_0$.
The proof of these results is presented in Section \ref{sec3}. Finally, in the
last section, we present a few consequences of the main result and we
detail some known and new examples. We also mention that some of the
results stated in Theorem~\ref{thm:2} below were announced without
proofs in the note \cite{Patie-09-cras}.

\section{Main results}\label{sec2}
Henceforth, we assume that $X$ is a pssMp of index $1/\alpha>0$ and of
the spectrally negative type. It means that it is associated via the
Lamperti mapping to a possibly killed L\'{e}vy process $\xi$ which is
spectrally negative. We exclude the cases when $\xi$ is degenerate,
that is, when $\xi$ is the negative of a subordinator or a pure drift
process. We recall that $\Pl$ (resp., $\El$) stands for the law
(resp., the expectation operator) of $\xi$ with $\xi_0=0$. The law of
$\xi$ is determined by its Laplace exponent $\overline{\psi}(u)=\psi
(u)-q$, where $q\geq0$ is the killing rate and $\psi$ admits the
following L\'{e}vy--Khintchine representation: for any $u\geq0$,
\[
\psi(u) = \bar{b} u + \frac{\sigma}{2} u^2 + \int_{-\infty}^0
\bigl(e^{u r} -1
-ur{{\mathbb{I}}}_{\{|r|<1\}} \bigr)\nu(dr),
\]
where $ \bar{b}\in\R, \sigma\geq0$ and the measure
$\nu$ is such that $\int_{-\infty}^0 (1
\wedge r^2 ) \nu(dr) <+ \infty$.
We shall refer to $\xi$ (resp., $\overline{\psi}$) as the underlying
L\'
{e}vy process (resp., Laplace exponent) of $X$.
Let us now proceed by recalling some basic properties of the Laplace
exponent $\psi$, which can be found, for instance, in Bertoin \cite
{Bertoin-96}. First, it is plain that $\lim_{u \rightarrow
\infty}\psi(u)=+\infty$ and by monotone convergence, one gets
$\El[\xi_1]= \bar{b} +\int_{-\infty}^{-1}r\nu(dr) \in[-\infty
,\infty)$.
We shall also need the value of the constant $\Lambda=\lim
_{u\rightarrow\infty} \frac{\psi(\alpha u)}{u}$ which is given (see
\cite{Bertoin-96}, Corollary VII.5) by
\[
\Lambda= \cases{
\displaystyle \alpha b =\alpha\biggl(\bar{b}-\int_{-1}^0 r \nu(dr)\biggr), &\quad if
$\sigma=0$ and $\displaystyle \int_{-\infty}^0 (1 \wedge
r)\nu(dr)<\infty$,\vspace*{2pt}\cr
+\infty, &\quad otherwise.}
\]
Since we have excluded the degenerate cases, we easily check that
$b>0$. Next, we recall that the\vadjust{\goodbreak} mapping $u\mapsto\psi(u)$ is
continuous and increasing on $[\phi(0),\infty)$, where $\phi(0)$ stands
for the largest solution to the equation $\psi(u)=0$. Thus, $\psi$
has a
well-defined inverse function $\phi\dvtx[0,\infty)\rightarrow
[\phi(0),\infty)$ which is also continuous and increasing. In order to
simplify the notation we write, for any $q\geq0$, $\gamma=\phi(q)>0$.
Then, it is easily seen that
\[
\El[e^{ \gamma\xi_1}] = 1.
\]
We also note that the condition \hyperlink{HypoH}{H} is equivalent to the
requirement $\phi(q)>0$.
Next, we set $\psi_{\gamma}(u)=\psi(u+\gamma)-\psi(\gamma)$ and
observing that $\psi_{\gamma}(0)=0$, we deduce that $\psi_{\gamma}$ is
the Laplace exponent of a conservative spectrally negative L\'{e}vy
process. We also point out that $\psi_{\gamma}'(0^+)=\psi'(\gamma)>0$
and $\lim_{u\rightarrow\infty} \frac{\psi_{\gamma}(u)}{u}=\lim
_{u\rightarrow\infty} \frac{\psi(u)}{u}$.

We proceed by introducing more notation taken from Patie
\cite{Patie-OU-06} and \cite{Patie-06c}. First, for a function $f$ and for
any $\alpha>0$, we write
\[
a_{s}(f;\alpha) = \prod_{k=1}^{\infty}\frac{f(\alpha(k
+s)
)}{f(\alpha k)}, \qquad  s\in\C,
\]
whenever the infinite product exists. Note that, for instance, $a_0(\psi
;\alpha)=1$ and for any $n=1,2,\ldots,$
%
\begin{equation} \label{eq:coef}
a_n(\psi;\alpha) =\Biggl( \prod_{k=1}^n\psi(\alpha k)\Biggr)^{-1}.
\end{equation}
Next, we introduce, for any $\rho\in\C$ such that $\mathfrak
{Re}(\rho)>0$, the
power series
%
\begin{equation} \label{eq:f1}
\Ip(\rho;z)= \frac{1}{\Gamma(\rho)}\sum_{n=0}^{\infty}
\ann
\Gamma(\rho+n) z^{n},
\end{equation}
where $\Gamma$ stands for the Gamma function. By means of classical
criteria, it is easily seen that the function $z \mapsto\Ip(\rho;z)$
is analytic in the disc $ \{z \in\C;\break |z|<\Lambda\}$. In particular,
in the case $\Lambda=+\infty$, that is, when the process $\xi$ has paths
of unbounded variations, $\Ip(\rho;z)$ is an entire function in $z$.
Moreover, for any $|z|<\Lambda$, the mapping $\rho\mapsto\Ip(\rho;z)$
is a meromorphic function defined for all complex numbers $\rho$ except
at the poles of the Gamma function, which are the points $\rho
=0,-1,\ldots$ However, they are removable singularities. Indeed, for
any $|z|<\Lambda$ and any integer $N\in\Nn$, one has, by means of the
recurrence relation $\Gamma(z+1)=z\Gamma(z)$,
\[
\Ip(0;z)=1
\]
and
\[
\Ip(-N;z) = \sum_{n=0}^{N}(-1)^n \frac{\Gamma(N+1)}{\Gamma
(N+1-n)}\ann z^{ n}.
\]
Thus, by uniqueness of the analytic continuation, for any $|z|<\Lambda
$, $\Ip(\rho;z)$ is an entire function in $\rho$. Before stating our
main result,\vadjust{\goodbreak} we show that in the case $\Lambda=\alpha b$, the power
series (\ref{eq:f1}) can be represented, in the left half-plane, as
another convergent power series which corresponds to an analytic
continuation in this domain. To this end, we aim to use the co-called
Euler transformation; see, for example, \cite{Norlund-55}, page 294.
However, this transformation can be performed if and only if the
singularity of the function $\Ipg(\rho;z)$ on the circle $|z|=\Lambda$
is located at the point $z=\Lambda$. In order to show that our family
of functions satisfies this property, we first provide a contour
integral representation of $\Ipg(\rho;z)$ which turns out to be an
analytic continuation in the entire complex plane cut along the
positive real axis. Then, we are able to apply the Euler
transformation to derive a series representation.\vspace*{-2pt}
\begin{prop} \label{prop:ac}
Let $\Lambda=\alpha b$, then $\Ipg(\rho;z)$ is analytic in the disc
$|z|<\alpha b$ and for any fixed $\rho=0,-1,\ldots,$ the mapping
$z\mapsto\Ipg(\rho;z)$, as a~polynomial, is an entire function.

Moreover, for any $\rho\neq0,-1,\ldots, \Ipg(\rho;z)$ admits an
analytic continuation in the entire complex plane cut along the
positive real axis given by
%
\begin{eqnarray} \label{eq:ci}
\Ipg(\rho;z) = \frac{1}{2i \pi\Gamma(\rho)}\int_{-i\infty
}^{i\infty}
a_s(\varphi_{\gamma};\alpha) \Gamma(s+\rho)\Gamma(-s)
\biggl(-\frac
{z}{\alpha}\biggr)^{s} \,ds,\nonumber\\[-8pt]\\[-8pt]
&&\eqntext{|{\arg}(-z)|<\pi,}
\end{eqnarray}
where the contour is indented to ensure that all poles
(resp., nonnegative poles) of $\Gamma(\rho+s)$ [resp., $\Gamma(-s)$] lie
to the left (resp., right) of the intended imaginary axis.

Consequently, for any $\rho\in\C$, $\Ipg(\rho;z)$ admits, in the
half-plane $\mathfrak{Re}(z)<\frac{\alpha b}{2}$, the following power series
representation
%
\begin{equation}\label{eq:h}\qquad
\Ipg(\rho;z) = \biggl(1-\frac{z}{\alpha b}\biggr)^{-\rho} \sum
_{n=0}^{\infty} \Ipg(-n;\alpha b)\frac{\Gamma(\rho+n)}{n!\Gamma
(\rho
)}\biggl(\frac{z}{z-\alpha b}\biggr)^n.
\end{equation}
Finally, for any fixed $\mathfrak{Re}(z)<\frac{\alpha b}{2}$, $\Ipg
(\rho;z)$ is
an entire function in the argument~$\rho$.\vspace*{-2pt}
\end{prop}
\begin{remark}
A specific instance of the mapping $\I_{\psi_{\gamma}}(\rho;x)$ when
$\Lambda=\alpha b$, is the hypergeometric function ${}_2F_1$. In this
case, the representation (\ref{eq:h}) is known as the Euler
transformation which has the remarkable feature that the power series
on the right-hand side of (\ref{eq:h}) is still an hypergeometric
function ${}_2F_1$. We refer to the Section \ref{ex:3} below for more
details on this example.\vspace*{-2pt}
\end{remark}

We are now ready to state our main result.\vspace*{-2pt}
\begin{theorem} \label{thm:2}
Let $q\geq0$, assume that $\phi(q)>0$ and set $\gamma=\phi(q)$ and
$\gamma_{\alpha}=\gamma/\alpha$.
Then, there exists a constant $C_{\gamma}>0$ such that
%
\begin{equation} \label{eq:const1}
\I_{\psi_{\gamma}}(\gamma_{\alpha};-t) \sim \frac{t^{-\gamma
_{\alpha
}}}{C_{\gamma}}  \qquad\mbox{as }   t\rightarrow\infty\vadjust{\goodbreak}
\end{equation}
($f(t)\sim g(t) $ as $t\rightarrow a$ means that $\lim_{t\rightarrow
a}\frac{f(t)}{g(t)}=1$ for any $a \in[0,\infty]$) and
%
\begin{equation}\label{eq:pd}
S(t) = C_{\gamma}t^{-\gamma_{\alpha}} \I_{\psi_{\gamma}}(\gamma
_{\alpha
};-t^{-1}),\qquad  t>0,
\end{equation}
where, by self-similarity, we have set $S(tx^{-\alpha})=\Q
_{x}(T_0\geq
t), x,t>0$.
Finally, the law of $T_0$ under $\Q_1$ is absolutely continuous with an
infinitely continuously differentiable density denoted by $s$ and given by
\[
s(t) = \gamma_{\alpha} C_{\gamma} t^{-\gamma_{\alpha}-1} \I
_{\psi
_{\gamma}}(1+\gamma_{\alpha};-t^{-1}), \qquad  t>0.
\]
\end{theorem}
\begin{remark}
In the case $\Lambda=\infty$, we easily check that, for any
$\mathfrak{Re}
(\rho
)>0$, the mapping $x\mapsto\I_{\psi_{\gamma}}(\rho;x)$ is increasing
on $[0,\infty)$. Hence, we deduce from the above theorem that the
entire function $z\mapsto\I_{\psi_{\gamma}}(\gamma_{\alpha};z)$
has no
real zeros.
\end{remark}

In the above theorem, the constant $C_{\gamma}$ is characterized by the
behavior of the function $\I_{\psi_{\gamma}}(\gamma_{\alpha};-t)$ for
large values of $t$. In what follows, we provide some representations
of this constant in terms of the Laplace exponent $\psi_{\gamma}$.
\begin{prop} \label{prop:ck}
\begin{longlist}[(1)]
\item[(1)] If $\Lambda=\alpha b$, then
\[
C_{\gamma} = \alpha^{\gamma_{\alpha}} a_{-\gamma_{\alpha
}}(\varphi
_{\gamma};\alpha),
\]
where $\varphi_{\gamma}(u)= b- \int_0^{\infty}e^{-ur}\int_{-\infty
}^{-r}e^{\gamma v}\nu(dv)\,dr$.
\item[(2)] Otherwise, we have
\[
C_{\gamma} =\cases{\psi'_{\gamma}(0^+), &\quad if $\gamma_{\alpha} =
1$,\cr
\displaystyle \alpha^n\psi'_{\gamma}(0^+)\Biggl(\prod_{k=1}^n\varphi_{\gamma
}(\alpha
k)\Biggr)^{-1}, &\quad if $\gamma_{\alpha} = n+1,  n=1,2\ldots,$
\cr
\displaystyle \frac{ \alpha^{2\gamma_{\alpha}}}{\Gamma(1-\gamma_{\alpha
})}a_{-\gamma
_{\alpha}}(\bar{\varphi}_{\gamma};\alpha), &\quad otherwise,}
\]
where $\varphi_{\gamma}(\alpha u)= \psi_{\gamma}(\alpha u) / \alpha
u$ and
\[
\bar{\varphi}_{\gamma}( u)=\frac{\hat{b}}{ u} + \frac{\sigma}{2}
+ \int
^{\infty}_0 e^{- u r}\int_{-\infty}^{-r}\int_{-\infty
}^{-s}e^{\gamma
v}\nu(dv) \,ds \,dr
\]
with $\hat{b}=\bar{b}+\sigma\gamma+\int_{-\infty}^0
(e^{\gamma
r}-{{\mathbb{I}}}_{\{|r|<1\}})r\nu(dr)$.
\item[(3)] Finally, if $q=0$ and $0<\gamma_{\alpha}<1$, then
\[
\Ip(r)\sim C_{\gamma}   \Gamma(1-\gamma_{\alpha} )
r^{\gamma_{\alpha}}\I_{\psi_{\gamma}}(r)
\qquad\mbox{as }   r\rightarrow
\infty,
\]
where $\Ip(r)= \sum_{n=0}^{\infty}
\ann r^{n}$ is an entire function.
\end{longlist}
\end{prop}

\section{Proofs}\label{sec3}
\subsection{A useful analytic continuation} \label{sec:pt1}
The first claim of Proposition \ref{prop:ac} follows from the
discussion preceding the proposition. Thus, let us assume that\vadjust{\goodbreak}
$\rho\neq0,-1,\ldots$ Since $\psi_{\gamma}'(0^+)>0$, $\psi
_{\gamma}$
is well defined and analytic in the positive right half-plane and $\psi
_{\gamma}(u)>0$ for any $u>0$. Our next aim is to extend the
coefficients $a_n(\psi_{\gamma},\alpha)$ to a function of the complex
variable. Since the paths of the L\'{e}vy process $\xi$ are of bounded
variation, its Laplace exponent $\psi_{\gamma}$ admits the following
representation (see \cite{Bertoin-96}, Section VII.3):
\[
\psi_{\gamma}(u)= u\bigl(b-\hat{v}_{\gamma}(u)\bigr),
\]
where $\hat{v}_{\gamma}(u)=\int_0^{\infty}e^{-ur}\int_{-\infty
}^{-r}e^{\gamma v}\nu(dv)\,dr$. Thus, for any $n\geq0$, we have
\[
a_n(\psi_{\gamma},\alpha)=\frac{1}{\Gamma(n+1)\alpha
^n}a_n(\varphi
_{\gamma};\alpha)
\]
with $a_n(\varphi_{\gamma};\alpha)^{-1}=\prod_{k=1}^n \varphi
_{\gamma
}(\alpha k)$ and $a_0(\varphi_{\gamma};\alpha)=1$.
It is plain that the mapping $\hat{v}_{\gamma}$ is analytic in
$F_{-\gamma}=\{s\in\C; \mathfrak{Re}(s)> -\gamma\}$ and $\hat
{v}_{\gamma
}(u)$ is decreasing on~$\R^+$ with $0<\hat{v}_{\gamma}(0)< b$ since
$\psi_{\gamma}(0^+)>0$. Then, we may write
\begin{eqnarray*}
a_{s}(\psi_{\gamma};\alpha)&=& \frac{1}{\Gamma(s+1)\alpha
^s}a_s(\varphi
_{\gamma};\alpha)\\
&=&\frac{1}{\Gamma(s+1)\alpha^s}\prod
_{k=1}^{\infty
}\frac{\varphi_{\gamma}(\alpha(k +s))}{\varphi_{\gamma}(\alpha k)},
\end{eqnarray*}
where the infinite product is easily seen to be absolutely convergent
for any $\mathfrak{Re}(s)>0$ by taking the logarithm and noting that
$|\hat
{v}_{\gamma}(s)|\leq\hat{v}_{\gamma}(\mathfrak{Re}(s))$; see, for
example, \cite{Titchmarsh-39}, Section 1.41. Moreover, $a_s(\varphi
_{\gamma};\alpha)$ satisfies the functional equation
\[
a_{s+1}(\varphi_{\gamma};\alpha)=\frac{1}{\varphi_{\gamma}(\alpha
(s+1))}a_{s}(\varphi_{\gamma};\alpha),
\]
which shows that $a_s(\varphi_{\gamma};\alpha)$ is analytic in the
half-plane $F_{-\gamma-1}=\{s\in\C;\break \mathfrak{Re}(s)> -1-\gamma\}$.
Consequently, $a_{s}(\varphi_{\gamma};\alpha)$ is bounded on any closed
subset of $F_{-\gamma-1}$.
Then, we set $G(s)=\Gamma(s+\rho)\Gamma(-s) a_{s}(\varphi_{\gamma
};\alpha)$ and define
\[
\mathfrak{I}_{\mathfrak{L}_R} =- \frac{1}{2i \pi\Gamma(\rho)}\int
_{\mathfrak{L}_R}G(s) \biggl(-\frac{z}{\alpha}\biggr)^s \,ds,
\]
where the integral is taken in a clockwise direction round the contour
$\mathfrak{L}_R$, consisting of a large semi-circle, of center the
origin and radius $R$, lying to the right of the imaginary axis. This
contour is intended to ensure that all poles (resp., nonnegative poles)
of $\Gamma(\rho+s)$ [resp., $\Gamma(- s)$] lie to the left (resp., right)
of the intended imaginary axis. This contour is always possible since
we have assumed that $\rho\neq0,-1,\ldots.$ We can split $\mathfrak
{I}_{\mathfrak{L}_R}$ up into two integrals,~$\mathfrak{I}_{\mathfrak
{A}_{iR}}$ along the imaginary axis and, writing $s=Re^{i\theta}$,
\[
\mathfrak{I}_{\mathfrak{C}_R}=-\frac{1}{2\pi i}\int_{-\pi/2}^{\pi
/2}G(Re^{i\theta}) \biggl(-\frac{z}{\alpha}\biggr)^{Re^{i\theta
}}Re^{i\theta} \,d\theta.
\]
Recalling the following well-known
asymptotic formulae (see,
e.g., \cite{Paris-Kaminski-01}, Section~2.4), as $|s|\rightarrow
\infty$,
\begin{eqnarray*}
\Gamma(s+ \rho) &\sim&\sqrt{2\pi}
e^{-Re^{i\theta}}R^{Re^{i\theta}+\rho-{1}/{2}}e^{i\theta
(Re^{i\theta
}+\rho-{1}/{2})}, \qquad |\theta| < \pi,\\
\Gamma(-s) &\sim&e^{-\pi R|{\sin\theta}|}
e^{Re^{i\theta}}R^{-Re^{i\theta}-{1}/{2}}e^{i\theta
(-Re^{i\theta
}-{1}/{2})},\qquad  |\theta| < \pi,
\end{eqnarray*}
and
\[
\biggl|\biggl(-\frac{z}{\alpha}\biggr)^{s}\biggr|\sim|\alpha
z|^{R\cos
\theta}e^{-R\sin\theta\arg(-z)},
\]
we deduce that as $|s|\rightarrow\infty$
%
\begin{eqnarray}\label{eq:ae}
\biggl|G(s)\biggl(-\frac{z}{\alpha}\biggr)^{s}\biggr| &\sim&
a R^{\mathfrak{Re}(\rho)-1}|\alpha z|^{R\cos\theta}\nonumber\\[-8pt]\\[-8pt]
&&{}\times\cases{ e^{-R
|{\sin
\theta}|
(\pi+\arg(-z))}, &\quad $0<\theta\leq
\pi/2$,\cr
e^{-R |{\sin\theta}| (\pi-\arg(-z))}, &\quad
$-\pi/2\leq\theta<0$,}\nonumber
\end{eqnarray}
where $a$ is a positive constant.
On the one hand, along the path ${\mathfrak{A}_{iR}}$ we have $\theta
=\pm\frac{\pi}{2}$ and thus as $|z|\rightarrow\infty$
\[
\biggl|G(s) \biggl(-\frac{z}{\alpha}\biggr)^{s}\biggr| \sim
a R^{\mathfrak{Re}(\rho)-1}e^{\pm({\pi}/{2}) \Im(\rho)} \cases
{ e^{-R
(\pi
+\arg(-z))}, &\quad
$\theta=\pi/2$,\cr
e^{-R (\pi-\arg(-z))}, &\quad $\theta=-\pi/2$.}
\]
For the integral (\ref{eq:ci}) to converge absolutely, it is therefore
required that $|{\arg}(-z)|<\pi$. On the other hand,
the asymptotic estimate (\ref{eq:ae}) gives, as $R\rightarrow\infty$,
\[
\mathfrak{I}_{\mathfrak{C}_R}\rightarrow0  \qquad\mbox{if } |z|<1
\mbox{ and } |{\arg}(-z)|<\pi.
\]
Thus, as $R\rightarrow\infty$,
\[
\mathfrak{I}_{\mathfrak{L}_R} \rightarrow-\frac{1}{2i \pi}\int
_{-i\infty}^{i\infty}G(s) \biggl(-\frac{z}{\alpha}\biggr)^{s} \,ds.
\]
Finally, evaluating $\mathfrak{I}_{\mathfrak{L}_R}$ by the Cauchy
integral theorem and letting $R\rightarrow\infty$, we get
\begin{eqnarray}
\frac{1}{2i\pi}\int_{-i\infty}^{i\infty}G(s) \biggl(-\frac
{z}{\alpha
}\biggr)^{s} \,ds=\frac{1}{\Gamma(\rho)}\sum_{n=0}^{\infty}
a_n(\psi_{\gamma},\alpha)
\Gamma(\rho+n) z^{n}, \nonumber\\
&&\eqntext{|z|<1 \mbox{ and } |{\arg}(-z)|<\pi.}
\end{eqnarray}
Therefore, the integral (\ref{eq:ci}) offers an analytic continuation
of the mapping $z\mapsto\Ipg(\rho;z)$ in the entire complex plane cut
along the positive real axis.
Moreover, we deduce from such an analytic continuation that the power
series (\ref{eq:f1}) has an unique singularity on the circle
$|z|=\alpha b$ located at the point $z=\alpha b>0$. Now, following a
device developed for hypergeometric series (see N{\o}rlund\vadjust{\goodbreak}
\cite{Norlund-55}, pages 294 and 295), we introduce the function $\mathcal{H}$
defined for some $a\in\C$ by
\[
\mathcal{H}_{\psi_{\gamma},a}(\rho;z) = (1-z)^{-\rho}
\Ipg
\biggl(\rho;\frac{a \alpha b z}{z-1} \biggr).\vspace*{-1pt}
\]
Note that
%
\begin{equation}\label{eq:hir}
\Ipg(\rho;az) = \biggl(1-\frac{z}{\alpha b}\biggr)^{-\rho}
\mathcal
{H}_{\psi_{\gamma},a}\biggl(\rho;\frac{z}{z-\alpha b} \biggr).\vspace*{-1pt}
\end{equation}
Thus, denoting by $(b_n)_{n\geq0}$ the coefficients of the power series
$\mathcal{H}_{\psi_{\gamma},a}(\rho;z)$, we have $b_0=a_0$ and by means
of residues calculus, with $\mathfrak{C}$ a circle around $0$ of small
radius and with positive orientation, we have for $n\geq1$,
\begin{eqnarray*}
b_n &=& \frac{1}{2\pi i } \int_{\mathfrak{C}}
\frac
{\mathcal{H}_{\psi_{\gamma},a}(\rho;z)}{z^{n+1}}\,dz
\\[-2pt]
&=& (-1)^n\frac{1}{2\pi i } \int_{\mathfrak{C}}(1-z
)^{-\rho} \Ipg\biggl(\rho;\frac{a \alpha b z}{z-1} \biggr) \frac
{dz}{z^{n+1}}
\,dv\\[-2pt]
&=& \frac{1}{\Gamma(\rho)}\sum_{k=0}^{n}(-a \alpha b)^{k}a_k(\psi
_{\gamma};\alpha)
\frac{\Gamma(\rho+n)}{\Gamma(n-k+1)}.\vspace*{-1pt}
\end{eqnarray*}
Thus, one gets
\[
(1-z)^{-\rho} \Ipg\biggl(\rho;\frac{a \alpha b z}{z-1} \biggr)
=\sum
_{n=0}^{\infty}\Ipg(-n;a \alpha b) \frac{\Gamma(\rho+n)}{\Gamma
(\rho
)n!} z^n.\vspace*{-1pt}
\]
From Weierstrass's double series theorem, the above identity is true if
$|z|<\frac{1}{1+|a|}$. Moreover, the function on the left-hand side has
a singularity at $z=1$ and $z=\frac{1}{1-a}$. Thus, the series on the
right-hand side is convergent if $|z|<1$ and $|z(1-a)|<1$. By choosing
$a=1$, we conclude by observing that the series on the right-hand side
of (\ref{eq:hir}) is convergent for $\mathfrak{Re}(z)<\frac{\alpha
b}{2}$.\vspace*{-2pt}

\subsection{The distribution of $T_0$} \label{sec:pt}
We proceed by introducing the Ornstein--Uhlenbeck process
$U=(U_t)_{t\geq
0}$ defined by
\[
U_t=e^{\tilde{\alpha}t}X_{\tau(t)},\qquad  t\geq0,\vspace*{-1pt}
\]
where $\ta=\alpha^{-1}$ and $\tau(t)=1-e^{-t}$. Next, we put
\[
H_0=\inf\{s>0;   U_s=0\}\vspace*{-1pt}
\]
and set
\[
1-K(x) = \Q_x(H_0<\infty), \qquad  x>0.\vspace*{-1pt}
\]
We are now ready to state the following.\vspace*{-2pt}
\begin{prop} \label{thm:1}
Assume that the condition \textup{\hyperlink{HypoH}{H}} holds. Then, for any $x>0$ and
$t>0$, we have
%
\begin{equation} \label{eq:ii}
K(xt^{-\tilde{\alpha}}) = \Q_x( T_0\geq t)\vspace*{-1pt}
\end{equation}
and $P$ is increasing on $\R^+$ with\vadjust{\goodbreak} $\lim_{x\rightarrow\infty
}K(x)=1$ and $K(0)=0$.
\end{prop}
\begin{pf}
First, a simple time change yields the following identity in
distribution:
\[
H_0\stackrel{(d)}{=}-\log(1-{T_0}\wedge1).
\]
Thus, we deduce that
\begin{eqnarray*}
1-K(x) &=&\Q_x(H_0<\infty)\\
&=&\Q_x( T_0<1).
\end{eqnarray*}
Then, invoking the self-similarity property of $X$ we obtain the identity
\[
\Q_x( T_0\geq t)=\Q_{xt^{-\tilde{\alpha}}}( T_0\geq1)
\]
from which we deduce the identity (\ref{eq:ii}) and the properties
stated on $P$.
\end{pf}

According to Proposition \ref{thm:1}, our goal now is to derive an
expression of the function $K(x)=1-\Q_{x}(H_0< \infty),   x>0$.
Relying on the following identity:
\[
K(x)= \lim_{a\rightarrow\infty} \lim_{q\rightarrow0} \E_{x}
\bigl[e^{-qH_a}\mathbb{I}_{\{H_a<H_0\}}\bigr],
\]
where $H_a = \inf\{s>0;  U_s\geq a\}$, the problem reduces to the
computation of the functional $\E_{x}[e^{-qH_a}\mathbb{I}_{\{
H_a<H_0\}}]$. Actually, for technical reasons, we must deal first
with the functional $\E^{(\gamma)}_{x}[e^{-qH_a}]$ which is
the Laplace transform of the first passage time above for the
Ornstein--Uhlenbeck process associated to the pssMp $X$ with underlying
Laplace exponent $\psi_{\gamma}$. Finally, by means of Doob h-transform
arguments, we will be able to relate the latter functional to the
former one.

We use the notation introduced in Theorem \ref{thm:2} and take first
$X$ with underlying Laplace exponent $\psi_{\gamma}$. We denote its law
(resp., its expectation operator) by $\Q^{(\gamma)}$ (resp., $\E
^{(\gamma
)}$). In order to simplify the notation we set, without loss of
generality, $\alpha=1$. We recall that $\psi_{\gamma}(0)=0$ and
$\psi
_{\gamma}(0^+)>0$ and hence the condition \hyperlink{HypoH}{H} does not hold.
Next, we simply write $Q^{(\gamma)}=(Q^{(\gamma)}_t)_{t\geq0}$ for the
semigroup of $X$, that is, for any bounded Borelian function~$g$ and~$t$, $x>0$, one has
\[
Q^{(\gamma)}_tg(x)=\E^{(\gamma)}_x[g(X_t)].
\]
From \cite{Bertoin-Yor-02-b}, we have that $Q^{(\gamma)}$ is a Feller
semigroup on $[0,\infty)$. Next, we say, for any $r\in\R$, that a
function $I$ is $r$-invariant for $Q^{(\gamma)}$ if
\[
e^{-rt}Q^{(\gamma)}_tI(x)=I(x),  \qquad x>0.
\]
We start with the following lemma which is obtained readily from
\cite{Patie-06c}, Theorem~1.
\begin{lemma} \label{lem:1}
For any $r>0$, the mapping $x\mapsto \Ipg(-rx)$ is $-r$-invariant for
$Q^{(\gamma)}$.\vadjust{\goodbreak}
\end{lemma}

Following a device developed by the author in \cite{Patie-08a}, we show
how to construct some specific time--space invariant functions for the
semigroup $Q^{(\gamma)}$ in terms of its $r$-invariant functions. We
now state the following result which is a slight generalization of
\cite{Patie-08a}, Theorem 1 and Corollary 3.2.
\begin{lemma} \label{prop:iou}
For any $\mathfrak{Re}(\rho)>0$, the mapping $x\mapsto\Ipg(\rho
;-x)$ satisfies
the identity, for any $0\leq t< 1$,
%
\begin{equation} \label{eq:ts}
(1- t)^{-\rho}Q^{(\gamma)}_t \bigl(d_{(1-t)^{-1}}\Ipg\bigr)
(\rho
;-x) = \Ipg(\rho;-x), \qquad  x>0,
\end{equation}
where $d_cf(x)=f(cx),  c>0$.
\end{lemma}

Next, we introduce the stopping time $ D_a$ defined, for any $a>0$, by
\[
D_a = \inf\{0<s\leq1;  X_s=a(1- s) \}.
\]
Writing $(a)_+=\max(a,0)$, we have
%
\begin{equation} \label{eq:hi}
e^{-H_a}\stackrel{(d)}{=}(1- D_a)_+
\end{equation}
and, in particular, for $a=0$, since $D_0\stackrel{(d)}{=}T_0\wedge1$,
we obtain
\[
e^{-H_0}\stackrel{(d)}{=}(1- T_0)_+.
\]
For any $a>0$, we set
\[
\kappa(a) = \inf\{\kappa\in\R^{+};  \Ipg(\kappa;-a)=0\}
\]
with the usual convention that $\inf\{\varnothing\}=\infty$, for the
smallest positive real zero of the function $\Ipg(\cdot;-a)$. We are now
ready to state the following.
\begin{cor} \label{cor:1}
Let $0\leq x\leq a$. Then, for any $\rho\in\C$ with $\mathfrak
{Re}(\rho
)<\kappa(a)$, we have
\[
\E^{(\gamma)}_x[(1 - D_a)_+^{-\rho}]=
\frac{\Ipg(\rho; -x)}{\Ipg(\rho;-a)}.
\]
Consequently, for any real $\kappa$ such that $\kappa<\km$, the mapping
$x\mapsto\Ipg(\kappa;-x)$ is positive on $\R^+$.
\end{cor}
\begin{pf}
Since $X$ under $\Q^{(\gamma)}$ is a Feller process on $[0,\infty)$,
we can start by fixing $x=0$ and $a>0$. Then, recalling that $\Ipg
(0,-a)=1$, we observe that $\Ipg(\kappa;-a )$ is positive for any
$0\leq\kappa< \km$ reals. The existence of such an interval follows
from the fact that the zeros of a nonconstant holomorphic function are
isolated. Thus, by combining the identity (\ref{eq:ts}) with the Dynkin
formula (see, e.g., \cite{Dynkin-65}, Theorem 12.4), applied to
the bounded stopping time~$D_a$, we deduce, for any $0\leq\kappa< \km
$, that
%
\begin{equation} \label{eq:ma}
\E^{(\gamma)}_0[(1- D_a)_+^{-\kappa}] = \frac{1}{\Ipg
(\kappa;-a )}.\vadjust{\goodbreak}
\end{equation}
Next, we recall, from   identity (\ref{eq:hi}), that
\[
e^{\kappa H_a}\stackrel{(d)}{=}(1- D_a)^{-\kappa}_+.
\]
Since $H_a$ is a positive random variable, as a Laplace transform, the
left-hand side on   identity (\ref{eq:ma}) is analytic in the
half-plane $\{\rho\in\C; \mathfrak{Re}(\rho)<\km\}$ and positive
on $\R
^+$; see,
for example, \cite{Widder-41}, Chapter II.
Then, let us assume that there exists a complex number $\rho(a)$ in the
strip $0\leq\mathfrak{Re}(\rho(a))<\km$ such that $\Ipg(\rho(a)
;-a )=0$.
However, as the left-hand side of (\ref{eq:ma}) is analytic with
respect to the argument $\kappa$ in this strip, we deduce, by the
principle of analytic continuation, that this is not possible.
Moreover, we get that $\Ipg(\rho;-a )$ has no zeros on $\{\rho\in\C;
 \mathfrak{Re}(\rho)<\km\}$ and is positive on $\{\kappa\in\R;\break
  \kappa
<\km\}
$. Finally, let us consider a real number $a_1$ such that $0< a_1\leq
a$. Clearly, $\Q^{(\gamma)}_0$-a.s. $(1- D_{a_1})_+^{-\kappa}\leq(1-
D_a)_+^{-\kappa}$, for any $0\leq\kappa< \km\wedge\kappa(a_1)$. Then
we deduce from~(\ref{eq:ma}), for any $0\leq\kappa< \km\wedge
\kappa(a_1)$, that
\[
0<\frac{1}{\Ipg(\kappa;-a_1 )} \leq\frac{1}{\Ipg(\kappa;-a )}.
\]
Thus, it is not difficult to see that $ \kappa(a_1) \geq\km$.
Therefore, since $ \kappa(x) \geq\km$, for any $0\leq x\leq a$, the
strong Markov property and the absence of positive jumps of $X$
complete the proof.
\end{pf}

The choice of starting our computation under the law $\Q^{(\gamma)}$
was motivated by the previous proof where it was necessary to start $X$
at $0$ in order to get some information about the sign of the function
$\Ipg(\kappa,-a)$. This device would not have been possible under $\Q$.
We proceed to the proof of Theorem~\ref{thm:2} which we now split into
two parts: the case when $X$ reaches $0$ continuously, that is, $q=0$
and $\El[\xi_1]<0$ and the case when $X$ reaches $0$ by a jump, that
is, $q>0$.

\subsubsection{Continuous killing}
Here, we assume that $q=0$ and $\El[\xi_1]<0$. Thus, in this case,
$\gamma=\phi(0)$ and $\psi_{\gamma}(u)=\psi(\gamma+u)$ with $\psi
_{\gamma}'(0^+)>0$.
\begin{lemma} \label{lem:2}
Writing $\kappa'(a)=\kappa(a)-\gamma>0$, we have, for any $\kappa
<\kappa'( a)$ and $0<x\leq a$,
\[
\E_x\bigl[(1- D_a)^{-\kappa}{{\mathbb{I}}}_{\{D_a<T_0\wedge1\}
}\bigr]=\frac
{x^{\gamma}}{a^{\gamma}}\frac{\Ipg(\kappa+ \gamma; -x )}{\Ipg
(\kappa
+\gamma;- a )}.
\]
In particular, for any $0<x\leq a$, we have
\[ 
\Q_x[ D_a<T_0\wedge1] = \frac{x^{\gamma}}{a^{\gamma
}}\frac
{\Ipg(\gamma; -x )}{\Ipg(\gamma;- a )}.
\]
\end{lemma}
\begin{pf}
We start by using the fact that the function $x\mapsto x^{-\gamma}$ is
excessive for $Q^{(\gamma)}_t$; see, for example, \cite{Rivero-05}. In
particular, one\vadjust{\goodbreak} has, for any $t>0$ and for any $F$ a $\F_t$-measurable
and bounded random variable,
\[
\E^{(\gamma)}_{x}[F] =\E_x[X_t^{\gamma} F, t<T_0
],\qquad  x>0.
\]
Note that this relation also holds for any $\F_{\infty}$-stopping time.
Moreover, proceeding as in the proof of Corollary \ref{cor:1}, one gets
that the Mellin transform of the positive random variable $(1- D_a)_+$
is well defined for any real $\kappa$ such that $\kappa\leq0$. Thus,
since $X$ has no positive jumps, one obtains by means of both Corollary
\ref{cor:1} and the optional stopping theorem, for any $\kappa\leq0$,
\begin{eqnarray*}
\E_x\bigl[(1 - D_a)_+^{-\kappa}{{\mathbb{I}}}_{\{D_a<T_0\}}
\bigr]&=&\frac
{x^{\gamma
}}{a^{\gamma}}\E^{(\gamma)}_x\bigl[(1 - D_a)_+^{-(\kappa+ \gamma
)}
\bigr] \\
&=& \frac{x^{\gamma}}{a^{\gamma}}\frac{\Ipg(\kappa+ \gamma; -x
)}{\Ipg
(\kappa+ \gamma; -a )}.
\end{eqnarray*}
We deduce that $\kappa'(a)>0$ and the proof is completed by letting
$\kappa\rightarrow0$.
\end{pf}

We are now ready to complete the proof of Theorem \ref{thm:2} in the
case $\gamma=\phi(0)$. One gets that
\begin{eqnarray*}
\Q_x[ D_a<T_0\wedge1]&=&\Q_x[\tau(H_a)<\tau
(H_0)\wedge
1]\\
&=& \Q_x[H_a<H_0]
\end{eqnarray*}
since $\tau$ is increasing and $\tau^{-1}(1)=\infty$. Thus, as $X$ has
no positive jumps, one deduces that
\begin{eqnarray*}
\lim_{a\rightarrow\infty}\Q_x[ D_a<T_0\wedge1]&=& \Q
_x
[H_0=\infty]\\
&=& K(x).
\end{eqnarray*}
As we have learnt from Corollary \ref{cor:1} and Lemma \ref{lem:2} that
the mapping $x\mapsto\Ipg( \gamma;-x )$ is positive on $\R^+$, it
means that there exists a constant $C_{\gamma}>0$ such that
\[
\Ipg( \gamma;-x ) \sim C_{\gamma}^{-1} x^{-\gamma} \qquad \mbox
{as }   x\rightarrow\infty.
\]
Then, recalling that $\lim_{x\rightarrow\infty}K(x)=1$, we obtain
\[
K(x)=C_{\gamma}x^{\gamma}\Ipg( \gamma; -x ).
\]
Hence, we deduce the expression of $S$ from the identity
$S(t)=K(t^{-1})$. Finally, the series $\Ipg( \gamma;-x )$ being
absolutely continuous, the expression of the density $s$ is obtained by
differentiating terms by terms. Indeed, one has
\begin{eqnarray*}
s(t)&=&-\frac{d}{dt}S(t)\\
&=& C_{\gamma}t^{-\gamma-1}\frac{1}{\Gamma( \gamma)}\sum
_{n=0}^{\infty}
(-1)^n a_n(\psi)(\gamma+n) \Gamma( \gamma+n)t^{-n}\\
&=&\frac{\Gamma(\gamma+1)}{\Gamma(\gamma)} C_{\gamma}t^{-\gamma
-1} \Ipg
(1+ \gamma;-t ).
\end{eqnarray*}
The expression of the successive derivatives are obtained by means of
an induction argument.

\subsubsection{$X$ reaches $0$ by a jump}
Throughout this part, we assume that $\xi$ is a spectrally negative L\'
{e}vy process killed at some independent exponential time of parameter
$q>0$. Recall that, for any $u\geq0$, $\tp(u)=\psi(u)-q$, $\phi$ is
such that $\psi\circ\phi(u)=u$ and with $\gamma=\phi(q)$, we
easily see
that $\psi_{\gamma}(u)=\tp(u+\gamma)$ and \mbox{$\psi'_{\gamma}(0^+)>0$}.
\begin{lemma}\label{lemmaa}
Writing $\kappa'(a)=\kappa(a)-\gamma>0$, we have, for any $\kappa
<\kappa'( a)$ and $0<x\leq a$,
\[
\E_x\bigl[(1- D_a)_+^{-\kappa}{{\mathbb{I}}}_{\{D_a<T_0\}}
\bigr]=\frac
{x^{\gamma
}}{a^{\gamma}}\frac{\Ipg(\kappa+ \gamma;- x )}{\Ipg(\kappa
+\gamma;- a )}.
\]
In particular,
\[
\Q_x[ D_a<T_0\wedge1 ]=\frac{x^{\gamma}}{a^{\gamma
}}\frac
{\Ipg(\gamma;-x )}{\Ipg(\gamma; -a )}.
\]
\end{lemma}
\begin{remark}\label{rem:1}
$\!\!$Writing $D_a^+=\inf\{s>0;  X_s=a(1+s)\}$ and $K(x;a)=\Q_x[
D_a<T_0\wedge1]$, we deduce from \cite{Patie-08a}, Corollary 3.2,
the following identity:
\[
\Q_x[ D^+_a<T_0]=K(-x,-a),  \qquad 0<x\leq a<\Lambda.
\]
It would be interesting to prove such a formula directly from the
definition of $D_a$ and $D_a^+$.
\end{remark}
\begin{pf*}{Proof of Lemma \ref{lemmaa}}
Let us observe from the Lamperti mapping~(\ref{eq:ss}) that the
semigroup $(Q_t)_{t\geq0}$ of $X$ is given for a function $f$
positive and measurable on $\R^+$ by
\[
Q_t f(x) = \E^{q}_x[e^{-q A_t}f(X_t)], \qquad  t\geq0,x>0,
\]
where $\E^{q}$ stands for the expectation operator associated to the
law of $X$ with underlying Laplace exponent $\psi$. Thus, for any $\F
_{\infty}$-stopping time $T$, one has
\[
\E_{x}[f(X_T)]=\E^{q}_x[e^{-q A_{T}}f(X_T)].
\]
Moreover, as $\xi$ has independent increments, it is plain that the
process $(e^{-qt+\gamma\xi_t})_{t\geq0}$ is a $\Pl^{q}$-martingale,
where $\Pl^{q}$ stands for the law of the L\'{e}vy process with Laplace
exponent $\psi$. By time change, one deduces that the process
$(X_t^{\gamma}e^{-qA_t})_{t\geq0}$ is a $\Q_1$-martingale. Thus, one
can define a new probability measure, which we denote by $\Q^{(\gamma
)}$, as follows, for any $t>0$ and for any $F$ a $\F_t$-measurable and
bounded random variable,
\[
\E^{(\gamma)}_x[F]=\E_x^{q}[X_t^{\gamma
}e^{-qA_t}F],\qquad  x>0.
\]
It is easily seen that the underlying Laplace exponent of $X$, under
$\Q
^{(\gamma)}$, is~$\psi_{\gamma}$.
Hence, one gets by the absence of positive jumps for $X$ and an
application of the optional stopping theorem, that, for any $0<x\leq a$
and $\kappa\leq0$,
\begin{eqnarray*}
\E_x\bigl[(1-D_a)_+^{-\kappa}{{\mathbb{I}}}_{\{D_a<T_0\}}
\bigr]&=& \E
^{q}_x
\bigl[e^{-q A_{D_a}}(1- D_a)_+^{-\kappa}{{\mathbb{I}}}_{\{D_a<T_0\}}
\bigr]\\
&=& \biggl(\frac{x}{a}\biggr)^{\gamma}\E^{(\gamma)}_x\bigl[(1-
D_a)_+^{-(\kappa+\gamma)}\bigr]\\
&=&\biggl(\frac{x}{a}\biggr)^{\gamma}\frac{\Ipg(\kappa+\gamma;-x
)}{\Ipg
(\kappa+\gamma;-a )},
\end{eqnarray*}
where the last line follows from Corollary \ref{cor:1} since $\psi
_{\gamma}'(0^+)>0$. The proof of the lemma is complete.
\end{pf*}

The proof of the theorem is completed by following a line of reasoning
similar to the previous case.

\subsection{\texorpdfstring{Proof of Proposition \protect\ref{prop:ck}}
{Proof of Proposition 2.5}}
Let us start by pointing out that it is not difficult to check that we
have, in all cases, $C_{\gamma}>0$. Moreover, let us first assume that
$\lim_{u\rightarrow\infty} \frac{\psi(u)}{u}=b$. From Proposition
\ref
{prop:ac}, we have
\begin{eqnarray}
\Ipg( \gamma_{\alpha} ;-z) =
\frac{1}{2i \pi\Gamma( \gamma_{\alpha} )}\int_{-i\infty
}^{i\infty}
a_s(\varphi_{\gamma};\alpha) \Gamma(s+ \gamma_{\alpha} )\Gamma
(-s)
\biggl(\frac{z}{\alpha}\biggr)^{s} \,ds, \nonumber\\
&&\eqntext{|{\arg}(z)|<\pi.}
\end{eqnarray}
Hence, upon displacement of the path to the left in order to include
the first pole of $\Gamma(s+ \gamma_{\alpha} )$ we obtain, from Theorem
\ref{thm:2} and a residue computation, that
\[
\Ipg( \gamma_{\alpha} ;-z) = \alpha^{ \gamma_{\alpha} }a_{-
\gamma
_{\alpha} }(\varphi_{\gamma};\alpha) z^{- \gamma_{\alpha} }+ o(z^{-
\gamma_{\alpha} }),
\]
which gives the characterization of $C_{\gamma}$ in this case.

For the other case, that is, when $\Lambda=+\infty$, one may follow a
line of reasoning similar to the proof of Proposition \ref{prop:ac}.
Indeed, as $0<\psi_{\gamma}'(0^+)<\infty$, we have, for any $u>0$,
\begin{eqnarray*}
\psi_{\gamma}(\alpha u) &=& \hat{b} \alpha u + \frac{\sigma}{2}
(\alpha
u)^2 + \int_{-\infty}^0 (e^{\alpha u r} -1
-\alpha ur)e^{\gamma r}\nu(dr) \\
&=& (\alpha u)^{2}\bar{\varphi}_{\gamma}(\alpha u),
\end{eqnarray*}
where $\hat{b}=\bar{b}+\sigma\gamma+\int_{-\infty}^0
(e^{\gamma
r}-{{\mathbb{I}}}_{\{|r|<1\}})r\nu(dr)$ and
\[
\bar{\varphi}_{\gamma}(\alpha u)=\frac{\hat{b}}{\alpha u} + \frac
{\sigma
}{2} + \int^{\infty}_0 e^{-\alpha u r}\int_{-\infty}^{-r}\int
_{-\infty
}^{-s}e^{\gamma v}\nu(dv) \,ds \,dr.
\]
Thus, as above, one may define the function
\begin{eqnarray*}
a_{s}(\psi_{\gamma};\alpha)&=& \frac{1}{\alpha^2\Gamma
^2(s+1)}a_s(\bar
{\varphi}_{\gamma};\alpha)\\[-2pt]
&=&\frac{1}{\alpha^2 \Gamma^2(s+1)}\prod_{k=1}^{\infty}\frac{\bar
{\varphi}_{\gamma}(\alpha(k +s+1))}{\bar{\varphi}_{\gamma}(\alpha k)}
\end{eqnarray*}
and observe the identity
\[
a_{s+1}(\bar{\varphi}_{\gamma};\alpha)=\frac{1}{\bar{\varphi
}_{\gamma
}(\alpha(s+1))}a_{s}(\bar{\varphi}_{\gamma};\alpha)
\]
with $a_0(\bar{\varphi}_{\gamma};\alpha)=1$. Hence, $a_{s}(\bar
{\varphi
}_{\gamma};\alpha)$ is a meromorphic function in $F_{-\gamma}=\{s\in
\C
;   \mathfrak{Re}(s)>-\gamma-1\}$ with simple poles at the points
$s_k=-k-1$ for
$k=0,1,\ldots$ and $s_k>-\gamma-1$. We obtain, writing $\overline
{G}(s)=\frac{a_s(\bar{\varphi}_{\gamma};\alpha)}{\Gamma(s+1)}
\Gamma(s+
\gamma_{\alpha} )\Gamma(-s)$, the following identity:
\[
\Ipg( \gamma_{\alpha} ;-z) = \frac{1}{2i \pi\Gamma( \gamma
_{\alpha}
)}\int_{-i\infty}^{i\infty} \overline{G}(s)\biggl(\frac{z}{\alpha
^2}\biggr)^{s} \,ds,
\]
which is now valid in the sector $ |{\arg}(z)|<\pi/2$. As above, after a
displacement of the path to the left in order to include the first pole
of $\Gamma(s+ \gamma_{\alpha} )$ we obtain, from Theorem \ref{thm:2}
and a residue computation, that
\begin{eqnarray*}
\Ipg( \gamma_{\alpha} ;-z) &=&
\frac{1}{\Gamma( \gamma_{\alpha} )}\Biggl(\sum_{k=1}^{[ \gamma
_{\alpha}
]}\frac{\Gamma(k)\operatorname{Res}_{j=-k}a_{j}(\bar{\varphi}_{\gamma
};\alpha
)}{\Gamma(1-k)} \biggl(\frac{z}{\alpha^2}\biggr)^{-k} \\[-2pt]
&&\hspace*{83.6pt}{}+
\operatorname{Res}_{s=- \gamma_{\alpha} }\overline{G}(s) \biggl(\frac{z}{\alpha
^2}\biggr)^{-s}\Biggr)\\[-2pt]
&&{}+ o(z^{- \gamma_{\alpha} }),
\end{eqnarray*}
where the sum is $0$ if $[ \gamma_{\alpha} ]$, the integer part of $
\gamma_{\alpha} $, is lower than $1$.
Since $a_{s}(\bar{\varphi}_{\gamma},\alpha)$ has a simple pole at
$j=-1,\ldots,-[ \gamma_{\alpha} ]$, the terms in the sum vanish. Hence,
if $ \gamma_{\alpha} $ is not an integer $\overline{G}(s)$ has a simple
pole at $- \gamma_{\alpha} $ and the expression of $C_{\gamma}$ follows
readily in this case. If $ \gamma_{\alpha} =n+1$, then $\overline
{G}(s)$ has a double pole at $-(n+1)$ and using the recurrence
relations of both the gamma function and $a_{s}(\bar{\varphi}_{\gamma
};\alpha)$, we deduce that
\begin{eqnarray*}
&&\operatorname{Res}_{s=-(n+1)}\overline{G}(s)\\[-2pt]
&&\qquad= \lim_{s\rightarrow
-n-1}\frac{d}{ds}\bigl( (s+n+1)^2 \overline{G}(s)\bigr)\\[-2pt]
&&\qquad= \lim_{s\rightarrow-n-1}\frac{d}{ds}\Biggl( \alpha^{-n-2} \prod
_{k=1}^{n}\varphi_{\gamma}\bigl(\alpha(s+k)\bigr)\psi_{\gamma}\bigl(\alpha
(s+n+1)\bigr)a_{s+n+1}(\bar{\varphi}_{\gamma};\alpha)\Biggr)\\[-2pt]
&&\qquad= \alpha^{-n-2}\Gamma(n+1)\psi'_{\gamma}(0^+)\prod
_{k=1}^{n}\varphi
(\alpha k)
\end{eqnarray*}
and the result follows.
The second part of the proposition is proved as follows. Let us recall
that in \cite{Patie-06c}, the expression of the Laplace transform of
$T_0$, in the case $\El[\xi_1]<0$, $q=0$ and $\gamma<\alpha$ is given
for any $r,x \geq0$ as follows:
%
\begin{equation} \label{eq:lt_exp}
\E_x[e^{-r T_0} ] = \N_{\psi,\gamma}(r x ),
\end{equation}
where
\[
\N_{\psi,\gamma}(r) =\Ip(r)-
C(\gamma)r^{ \gamma_{\alpha}}\It(r)
\]
and the positive constant $C(\gamma)$ is characterized by
\[
\Ip(r)\sim C(\gamma)
r^{\gamma_{\alpha}}\I_{\psi_{\gamma}}(r)
\qquad\mbox{as }   r\rightarrow
\infty.
\]
Next, let us write $\hat{F}(r)=\E_1[1-e^{-r T_0} ]$. Then,
from (\ref{eq:lt_exp}), one deduces easily that
\[
\hat{F}(r) \sim C(\gamma)r^{\gamma_{\alpha}}  \qquad\mbox{as }
r\rightarrow0,
\]
which is equivalent, according to Bingham, Goldie and Teugels \cite{Bingham-Goldie-Teugels-89},
Corollary~8.1.7, to
\[
S(t)\sim\frac{ C(\gamma)}{\Gamma(1-\gamma_{\alpha})} t^{-\gamma
_{\alpha
}}  \qquad\mbox{as }   t\rightarrow\infty,
\]
which completes the proof.

\section{Some final remarks and illustrative examples} \label{sec:exa}
We start by offering a~few consequences of Theorem \ref{thm:2}.
\begin{cor}\label{cor:dd}
With the notation used and introduced in Theorem~\ref{thm:2}, we have,
writing $s^{(m)}=\frac{d^m }{dt^m}s$,
\[
s^{(m)}(t)= (-1)^m \frac{\Gamma(m+1+ \gamma_{\alpha} )}{\Gamma(
\gamma
_{\alpha} )} C_{\gamma}t^{- \gamma_{\alpha} -1-m} \Ipg(m+1+ \gamma
_{\alpha} ;-t^{-1}), \qquad  t>0.
\]
Moreover,
\[
S(t)\sim C_{\gamma}t^{- \gamma_{\alpha} }  \qquad\mbox{as }
t\rightarrow\infty
\]
and, for any $m=0,1\ldots,$
\[ 
s^{(m)}(t) \sim (-1)^m C_{\gamma}\frac{\Gamma(m+1+ \gamma_{\alpha}
)}{\Gamma( \gamma_{\alpha} )} t^{- \gamma_{\alpha} -1-m}
\qquad\mbox{as }   t\rightarrow\infty.
\]
\end{cor}

As pointed out by several authors (see Carmona, Petit and Yor~\cite
{Carmona-Petit-Yor-97}, Rive\-ro~\cite{Rivero-05} and Maulik and Zwart
\cite{Maulik2006}) the study of the exponential functional is also
motivated by its connection to some interesting random affine equations
which have been deeply studied by Kesten \cite{Kesten-73}. Relying on a
result of Kesten, Rivero (\cite{Rivero-05}, Lemma 4) shows that there
exists a constant $C>0$ such that one has the following asymptotic behavior
\[
S(t) \sim Ct^{-\ta\gamma}  \qquad\mbox{as }   t\rightarrow\infty,
\]
whenever the L\'{e}vy process satisfies a set of conditions. As we have
excluded the case when $-\xi$ is a subordinator, it is not difficult to
verify that the L\'{e}vy processes we consider in this paper satisfy
Rivero's conditions. Hence, Theorem \ref{thm:2} and Proposition \ref
{prop:ck} offers several characterizations of the Kesten's constant. We
also point out that the asymptotic behavior of the density in Corollary
\ref{cor:dd} could not be deduced directly from Rivero's result since
we do not know whether or not the density is ultimately monotone.

\subsection{The Bessel processes} \label{sub:bess}
We consider $\xi$ to be a $2$-scaled Brownian motion with drift $2b
\in
\R$ and killed at some independent exponential time of parameter $q>0$,
that is, $\tp(u)=2u^2+2b u-q$ and $2\phi(q)=\sqrt{2q+b^2}-b$. Note that
$\psi_{\phi(q)}(u)=2u^2+(2b+\phi(q)) u$. Its associated
self-similar process~$X$ is well known to be a Bessel process of index
$b$ killed at a rate $q\int_0^{t}X_s^{-2}\,ds$. Moreover, we obtain,
setting $\varrho=b+2\phi(q)$,
\begin{eqnarray*}
\I_{\psi_{\phi(q)}}(\rho;-x)&=&\frac{\Gamma(\varrho
+1)}{\Gamma(\rho)}\sum_{n=0}^{\infty} (-1)^n\frac{\Gamma(\rho
+n)}{n!\Gamma(n+\varrho+1)}(x/2)^n\\
&=&\Phi(\rho, \varrho+1;-x/2),
\end{eqnarray*}
where $\Phi$ stands for the confluent hypergeometric function. We refer
to Lebedev (\cite{Lebedev-72}, Section 9) for useful properties of this
function. Next, using the following asymptotic:
\[
\Phi(\rho, \varrho+1;-x) \sim\frac{\Gamma(\varrho
+1)}{\Gamma
(\varrho+1-\rho)} x^{-\rho}  \qquad\mbox{as }   x\rightarrow
\infty,
\]
we get that $C_{\phi(q)}=\frac{\Gamma(\varrho+1-\phi(q))}{2^{\phi
(q)}\Gamma(\varrho+1)}$.
Thus, we obtain, recalling that, for any $q>0$, $\varrho-\phi
(q)=b+\phi(q)>0$,
\begin{eqnarray*}
s_{\phi(q)}(t)&=&\phi(q) \frac{\Gamma(\varrho+1-\phi(q))}{2^{\phi
(q)}\Gamma(\varrho+1)}t^{-\phi(q)-1}\Phi\bigl(1+\phi(q), \varrho
+1;-(2t)^{-1}\bigr)\\
&=&\frac{b+\phi(q)}{2^{\phi(q)}\Gamma(\phi(q))}t^{-\phi(q)-1}
\int
_0^{1}e^{-{u}/({2t})}(1-u)^{\varrho-\phi(q)-1}u^{\phi(q)}\,du,
\end{eqnarray*}
which is   expression  (5.a) in \cite{Yor-01}, page 105.
Considering now the case $q=0$ and $b<0$, we obtain readily that $\phi
(0)=-b$ and
\begin{eqnarray*}
s_{\phi(0)}(t)&=&\frac{2^b}{\Gamma(-b)}t^{b-1}\Phi\bigl(1-b,
1-b;-(2t)^{-1}\bigr)\\
&=&\frac{2^b}{\Gamma(-b)}t^{b-1}e^{-{1}/({2t})}.
\end{eqnarray*}
Hence, we deduce the well-known identity $(T_0,\Q_1)\stackrel{(d)}{=}
\frac{1}{2G_{-b}}$ where we recall that $G_{-b}$ stands for a Gamma
random variable of parameter $-b>0$.

\subsection{Law of the maximum of spectrally positive stable L\'{e}vy
processes} \label{sec:ms}
Let~$Z$ be an $\alpha$-stable spectrally negative L\'{e}vy process,
with $1<\alpha<2$. Let us denote by~$X$ the process $Z$ killed upon
entering\vspace*{1pt} into the negative half-line.~$X$ is then a pssMp. Next, we
denote by $\hat{Z}$ the dual of $Z$, that is, $\hat{Z}=-Z$ which is a
$\alpha$-stable spectrally positive L\'{e}vy process. Then, by means of
the translation invariance of L\'{e}vy processes, we deduce readily the
following identities:
\begin{eqnarray*}
\Q_x(T_0 \leq t)&=& \mathbb{P}_x\Bigl(\inf_{0<s\leq t}
Z_s\leq
0\Bigr)\\
&=& \mathbb{P}\Bigl(\max_{0<s\leq t} \hat{Z}_s\geq x \Bigr),
\end{eqnarray*}
which can be written as follows:
%
\begin{equation} \label{eq:sa}
\mathbb{P}\Bigl(\max_{0\leq s \leq t}\hat{Z_s}\geq x
\Bigr)=K(xt^{-\alpha
}), \qquad x,t>0.
\end{equation}
The Laplace exponent of the underlying L\'{e}vy process of $X$ has been
computed Patie \cite{Patie-CBI-09} in terms of the Pochhammer symbol.
Instead of using this expression, we follow an alternative route.
Indeed, in \cite{Patie-OU-06}, the author computed the unique
increasing invariant function, say $P_{+}$, of the Ornstein--Uhlenbeck
process defined by
\[
\tilde{U}_t=e^{- t/\alpha}X_{e^{t}-1}, \qquad  t\geq0.
\]
The function $P_+$, is given, with $C$ a constant to be determined and
writing $\ta=1/\alpha$, by
\begin{eqnarray*} \label{eq:lin}
P_+(x) &=& C x^{\alpha-1}\sum_{n=0}^{\infty} \frac{\Gamma
(n+1-\tilde
{\alpha})}{\Gamma(\alpha n +\alpha)}\alpha^n x^{\alpha n} \\&=& C
x^{\alpha-1}{}_2\Psi_1
\left(\matrix{
(1,1),(1,1-\ta) \cr
(\alpha,\alpha)}\bigg|
\alpha x^{\alpha}\right)
, \qquad  x \geq0,
\end{eqnarray*}
where ${}_2\Psi_1$ stands for the Wright hypergeometric function. From
Remark~\ref{rem:1}, we have $K(x)= K_{+}(e^{i\pi/\alpha} x)$. Note that
$K(0)=0$ and using the large asymptotic of the function ${}_2\Psi_1$
(details can be found in \cite{Patie-08-Asym}), we get as $x\rightarrow
\infty$,
\[
{}_2\Psi_1
\left(\matrix{(1,1),(1,1-\ta) \cr
(\alpha,\alpha)}
\bigg|-x^{\alpha}\right) \sim\biggl(\frac{\sin( \tilde{\alpha}\pi
)}{\pi
}\biggr)^{-1} x^{1-\alpha}.
\]
Hence, by setting $C=\frac{\sin( \tilde{\alpha}\pi)}{\pi}$, we obtain
the required condition\break $\lim_{x\rightarrow\infty} K(\infty)=1$ and
\[
K(x)=\frac{\sin( \tilde{\alpha}\pi)}{\pi}x^{\alpha-1}{}_2\Psi
_1\left(\matrix{
(1,1),(1,1-\tilde{\alpha}) \cr
(\alpha,\alpha)}\bigg|
-x^{\alpha}\right).
\]
Next, from   identity (\ref{eq:sa}), we find that
\[
\mathbb{P}\Bigl(\max_{0\leq s\leq1}\hat{Z}_s\geq x\Bigr)= P(x),
\]
where $\hat{Z} $ is a spectrally positive stable process of index
$\alpha$. Thus, by differentiating, one gets the following expression
for the density:
\[
k(x)
=\frac{\sin( \tilde{\alpha}\pi)}{\pi}x^{\alpha-2}{}_2\Psi
_1
\left(\matrix{
(1,1),(1,1-\tilde{\alpha}) \cr
(\alpha,\alpha-1)}
\bigg|  -x^{\alpha}\right),
\]
which is the expression found by Bernyk, Dalang and Peskir \cite{Bernyk-Dalang-Peskir-08},
Theorem 1.

\subsection{The self-similar saw-tooth processes} \label{ex:3}
Finally, we consider the so-called saw-tooth process introduced and
deeply studied by Carmona, Petit and Yor \cite{Carmona-Petit-Yor-98}. It is a
self-similar positive Markov process of index $\alpha=1$ with
underlying L\'{e}vy process the sum of a drift of parameter $b=1$ and
the negative of a compound Poisson process of parameter
$\beta>0$ whose jumps are exponentially distributed with parameter
$\delta+\beta-1>0$, that is,
\[
\psi(u) =u\frac{u+\delta-1 }{u+\delta+ \beta-1}, \qquad  u\geq0.
\]
Moreover, in \cite{Carmona-Petit-Yor-98}, the authors show that
\[
\phi(q) =\tfrac{1}{2}\bigl(q-(\delta-1)+\bar{\phi}(q)
\bigr), \qquad q\geq0,
\]
where $\bar{\phi}(q)=\sqrt{(q-(\delta-1))^{2}+4(\delta+\beta-1) q}$.
Let us proceed with the case \mbox{$q=0$}. Note, for $1-\beta<\delta<1$, that
$\gamma=1-\delta$ and
\[
\psi_{1-\delta}(u)=u\frac{u+1-\delta}{u+\beta}.
\]
Thus,
\[
a_n(\psi_{1-\delta},1)= \frac{\Gamma(n+1+\beta)\Gamma(2-\delta
)}{\Gamma
(1+\beta)\Gamma(n+1)\Gamma(n+2-\delta)},\qquad
a_0=1,
\]
and for $|z|<1$
\begin{eqnarray*}
\I_{\psi_{1-\delta}}(\rho;-z)&=&\frac{\Gamma
(2-\delta
)}{\Gamma(\rho)\Gamma(1+\beta)}\sum_{n=0}^{\infty} (-1)^n\frac
{\Gamma
(\rho+n)\Gamma(n+1+\beta)}{\Gamma(n+2-\delta)n!} z^n\\
&=&{}_2F_1(\rho,1+\beta,2-\delta;-z),
\end{eqnarray*}
where ${}_2F_1(a,b;x)$ stands for the hypergeometric function; see
Lebedev \cite{Lebedev-72}, Section 9, for a detailed account on this
function. Next, recalling the identity
\[
{}_2F_1(-n,1+\beta,\delta;1) = \frac{\Gamma(2-\delta)\Gamma
(n+1-\delta
-\beta)}{\Gamma(2-\delta+n)\Gamma(1-\delta-\beta)},
\]
we recover from (\ref{eq:h}) the well-known identity
\[
{}_2F_1(\rho,1+\beta,2-\delta;z) = (1-z)^{-\rho}{}_2F_1\biggl(\rho
,1-\delta-\beta,\delta;\frac{z}{z-1}\biggr),
\]
which provides an analytic continuation of the hypergeometric function
into the half-plane $\mathfrak{Re}(z)<\frac{1}{2}$. Finally,
using the  asymptotic
\[
{}_2F_1(\rho,1+\beta,2-\delta;-x) \sim\frac{\Gamma(2-\delta
)\Gamma
(1+\beta-\rho)}{\Gamma(2-\delta-\rho)\Gamma(1+\beta)} x^{-\rho}
\qquad \mbox{as }   x\rightarrow\infty,
\]
one obtains
\[
S(t)=\frac{\Gamma(1+\beta)}{\Gamma(2-\delta)\Gamma(\beta
+\delta)}
t^{\delta-1}{}_2F_1(1-\delta,1+\beta,\delta;-t^{-1}).
\]
Moreover, after some easy computations, one gets for
$\gamma=\phi(q),q>0$,
\[
\psi_{\phi(q)}(u)=u\frac{u+\phi(q)}{u+\beta+\delta+\phi(q)-1}.
\]
Thus, proceeding as above, we obtain
\[
\I_{\psi_{\phi(q)}}(\rho;-z)={}_2F_1\bigl(\rho
,\beta+\delta
+\phi(q),1+\phi(q);-z\bigr)
\]
and
\begin{eqnarray*}
S(t)&=&\frac{\Gamma(\beta+\delta+\phi(q))\Gamma(1+\bar{\phi
}(q)-\phi
(q))}{\Gamma(1+\bar{\phi}(q))\Gamma(\beta+\delta)}\\
&&{}\times t^{-\phi
(q)}{}_2F_1\bigl(\phi(q),\beta+\delta+\phi(q),1+\bar{\phi
}(q);-t^{-1}\bigr).
\end{eqnarray*}

\section*{Acknowledgment}
I am grateful to M. Savov and an anonymous referee for their comments which
significantly helped in improving the presentation of the paper.

%

%
\printaddresses


\begin{thebibliography}{32}

\bibitem{Alberts-Sheffield-08}
%
\begin{barticle}[mr]
\bauthor{\bsnm{Alberts},~\bfnm{T.}\binits{T.}} \AND
\bauthor{\bsnm{Sheffield},~\bfnm{S.}\binits{S.}}
(\byear{2011}).
\btitle{The covariant measure of SLE on the boundary}.
\bjournal{Probab. Theory Related Fields}
\bvolume{149}
\bpages{331--371}.
\bid{mr={2776618}}
\end{barticle}
%
\endbibitem

\bibitem{Bernyk-Dalang-Peskir-08}
%
\begin{barticle}[mr]
\bauthor{\bsnm{Bernyk},~\bfnm{Violetta}\binits{V.}},
\bauthor{\bsnm{Dalang},~\bfnm{Robert~C.}\binits{R.~C.}} \AND
\bauthor{\bsnm{Peskir},~\bfnm{Goran}\binits{G.}}
(\byear{2008}).
\btitle{The law of the supremum of a stable {L}\'evy process with no negative
jumps}.
\bjournal{Ann. Probab.}
\bvolume{36}
\bpages{1777--1789}.
\bid{doi={10.1214/07-AOP376}, mr={2440923}}
\end{barticle}
%
\endbibitem

\bibitem{Bertoin-96}
%
\begin{bbook}[mr]
\bauthor{\bsnm{Bertoin},~\bfnm{Jean}\binits{J.}}
(\byear{1996}).
\btitle{L\'evy Processes}.
\bseries{Cambridge Tracts in Mathematics}
\bvolume{121}.
\bpublisher{Cambridge Univ. Press}, \baddress{Cambridge}.
\bid{mr={1406564}}
\end{bbook}
%
\endbibitem

\bibitem{Bertoin-02-f}
%
\begin{barticle}[mr]
\bauthor{\bsnm{Bertoin},~\bfnm{Jean}\binits{J.}}
(\byear{2002}).
\btitle{Self-similar fragmentations}.
\bjournal{Ann. Inst. H. Poincar\'e Probab. Statist.}
\bvolume{38}
\bpages{319--340}.
\bid{doi={10.1016/S0246-0203(00)01073-6}, mr={1899456}}
\end{barticle}
%
\endbibitem

\bibitem{Bertoin-Yor-02-b}
%
\begin{barticle}[mr]
\bauthor{\bsnm{Bertoin},~\bfnm{Jean}\binits{J.}} \AND
\bauthor{\bsnm{Yor},~\bfnm{Marc}\binits{M.}}
(\byear{2002}).
\btitle{The entrance laws of self-similar {M}arkov processes and exponential
functionals of {L}\'evy processes}.
\bjournal{Potential Anal.}
\bvolume{17}
\bpages{389--400}.
\bid{doi={10.1023/A:1016377720516}, mr={1918243}}
\end{barticle}
%
\endbibitem

\bibitem{Bertoin-Yor-05}
%
\begin{barticle}[mr]
\bauthor{\bsnm{Bertoin},~\bfnm{Jean}\binits{J.}} \AND
\bauthor{\bsnm{Yor},~\bfnm{Marc}\binits{M.}}
(\byear{2005}).
\btitle{Exponential functionals of {L}\'evy processes}.
\bjournal{Probab. Surv.}
\bvolume{2}
\bpages{191--212 (electronic)}.
\bid{doi={10.1214/154957805100000122}, mr={2178044}}
\end{barticle}\vadjust{\goodbreak}
%
\endbibitem

\bibitem{Bingham-Goldie-Teugels-89}
%
\begin{bbook}[mr]
\bauthor{\bsnm{Bingham},~\bfnm{N.~H.}\binits{N.~H.}},
\bauthor{\bsnm{Goldie},~\bfnm{C.~M.}\binits{C.~M.}} \AND
\bauthor{\bsnm{Teugels},~\bfnm{J.~L.}\binits{J.~L.}}
(\byear{1989}).
\btitle{Regular Variation}.
\bseries{Encyclopedia of Mathematics and Its Applications}
\bvolume{27}.
\bpublisher{Cambridge Univ. Press}, \baddress{Cambridge}.
\bid{mr={1015093}}
\end{bbook}
%
\endbibitem

\bibitem{Carmona-Petit-Yor-97}
%
\begin{bincollection}[mr]
\bauthor{\bsnm{Carmona},~\bfnm{Philippe}\binits{P.}},
\bauthor{\bsnm{Petit},~\bfnm{Fr{\'e}d{\'e}rique}\binits{F.}} \AND
\bauthor{\bsnm{Yor},~\bfnm{Marc}\binits{M.}}
(\byear{1997}).
\btitle{On the distribution and asymptotic results for exponential functionals
of {L}\'evy processes}.
In \bbooktitle{Exponential Functionals and Principal Values Related to
{B}rownian Motion}.
\bseries{Bibl. Rev. Mat. Iberoamericana}
\bpages{73--130}.
\bpublisher{Rev. Mat. Iberoamericana}, \baddress{Madrid}.
\bid{mr={1648657}}
\end{bincollection}
%
\endbibitem

\bibitem{Carmona-Petit-Yor-98}
%
\begin{barticle}[mr]
\bauthor{\bsnm{Carmona},~\bfnm{Philippe}\binits{P.}},
\bauthor{\bsnm{Petit},~\bfnm{Fr{\'e}d{\'e}rique}\binits{F.}} \AND
\bauthor{\bsnm{Yor},~\bfnm{Marc}\binits{M.}}
(\byear{1998}).
\btitle{Beta-gamma random variables and intertwining relations between certain
{M}arkov processes}.
\bjournal{Rev. Mat. Iberoamericana}
\bvolume{14}
\bpages{311--367}.
\bid{mr={1654531}}
\end{barticle}
%
\endbibitem

\bibitem{Dufresne-90}
%
\begin{barticle}[mr]
\bauthor{\bsnm{Dufresne},~\bfnm{Daniel}\binits{D.}}
(\byear{1990}).
\btitle{The distribution of a perpetuity, with applications to risk
theory and
pension funding}.
\bjournal{Scand. Actuar. J.}
\banumber{no. 1--2},
\bpages{39--79}.
\bid{mr={1129194}}
\end{barticle}
%
\endbibitem

\bibitem{Dynkin-65}
%
\begin{bbook}[mr]
\bauthor{\bsnm{Dynkin},~\bfnm{E.~B.}\binits{E.~B.}}
(\byear{1965}).
\btitle{Markov Processes. {V}ols. {I}, {II}}.
\bseries{Translated with the Authorization and Assistance of the
Author by J.
Fabius, V. Greenberg, A. Maitra, G.~Majone. Die Grundlehren der
Mathematischen Wissenschaften, B\"ande 121}
\bvolume{122}.
\bpublisher{Academic Press}, \baddress{New York}.
\bid{mr={0193671}}
\end{bbook}
%
\endbibitem

\bibitem{Gjessing-Paulsen-97}
%
\begin{barticle}[mr]
\bauthor{\bsnm{Gjessing},~\bfnm{H{\aa}kon~K.}\binits{H.~K.}} \AND
\bauthor{\bsnm{Paulsen},~\bfnm{Jostein}\binits{J.}}
(\byear{1997}).
\btitle{Present value distributions with applications to ruin theory and
stochastic equations}.
\bjournal{Stochastic Process. Appl.}
\bvolume{71}
\bpages{123--144}.
\bid{doi={10.1016/S0304-4149(97)00072-0}, mr={1480643}}
\end{barticle}
%
\endbibitem

\bibitem{Kesten-73}
%
\begin{barticle}[mr]
\bauthor{\bsnm{Kesten},~\bfnm{Harry}\binits{H.}}
(\byear{1973}).
\btitle{Random difference equations and renewal theory for products of random
matrices}.
\bjournal{Acta Math.}
\bvolume{131}
\bpages{207--248}.
\bid{mr={0440724}}
\end{barticle}
%
\endbibitem

\bibitem{Lamperti-72}
%
\begin{barticle}[mr]
\bauthor{\bsnm{Lamperti},~\bfnm{John}\binits{J.}}
(\byear{1972}).
\btitle{Semi-stable {M}arkov processes. {I}}.
\bjournal{Z. Wahrsch. Verw. Gebiete}
\bvolume{22}
\bpages{205--225}.
\bid{mr={0307358}}
\end{barticle}
%
\endbibitem

\bibitem{Lebedev-72}
%
\begin{bbook}[mr]
\bauthor{\bsnm{Lebedev},~\bfnm{N.~N.}\binits{N.~N.}}
(\byear{1972}).
\btitle{Special Functions and Their Applications}.
\bpublisher{Dover}, \baddress{New York}.
\bid{mr={0350075}}
\end{bbook}
%
\endbibitem

\bibitem{Matsumoto-Yor-05-1}
%
\begin{barticle}[mr]
\bauthor{\bsnm{Matsumoto},~\bfnm{Hiroyuki}\binits{H.}} \AND
\bauthor{\bsnm{Yor},~\bfnm{Marc}\binits{M.}}
(\byear{2005}).
\btitle{Exponential functionals of {B}rownian motion. {I}.
{P}robability laws
at fixed time}.
\bjournal{Probab. Surv.}
\bvolume{2}
\bpages{312--347 (electronic)}.
\bid{doi={10.1214/154957805100000159}, mr={2203675}}
\end{barticle}
%
\endbibitem

\bibitem{Matsumoto-Yor-05-2}
%
\begin{barticle}[mr]
\bauthor{\bsnm{Matsumoto},~\bfnm{Hiroyuki}\binits{H.}} \AND
\bauthor{\bsnm{Yor},~\bfnm{Marc}\binits{M.}}
(\byear{2005}).
\btitle{Exponential functionals of {B}rownian motion. {II}. {S}ome related
diffusion processes}.
\bjournal{Probab. Surv.}
\bvolume{2}
\bpages{348--384 (electronic)}.
\bid{doi={10.1214/154957805100000168}, mr={2203676}}
\end{barticle}
%
\endbibitem

\bibitem{Maulik2006}
%
\begin{barticle}[mr]
\bauthor{\bsnm{Maulik},~\bfnm{Krishanu}\binits{K.}} \AND
\bauthor{\bsnm{Zwart},~\bfnm{Bert}\binits{B.}}
(\byear{2006}).
\btitle{Tail asymptotics for exponential functionals of {L}\'evy processes}.
\bjournal{Stochastic Process. Appl.}
\bvolume{116}
\bpages{156--177}.
\bid{doi={10.1016/j.spa.2005.09.002}, mr={2197972}}
\end{barticle}
%
\endbibitem

\bibitem{Norlund-55}
%
\begin{barticle}[mr]
\bauthor{\bsnm{N{\o}rlund},~\bfnm{N.~E.}\binits{N.~E.}}
(\byear{1955}).
\btitle{Hypergeometric functions}.
\bjournal{Acta Math.}
\bvolume{94}
\bpages{289--349}.
\bid{mr={0074585}}
\end{barticle}
%
\endbibitem

\bibitem{Paris-Kaminski-01}
%
\begin{bbook}[mr]
\bauthor{\bsnm{Paris},~\bfnm{R.~B.}\binits{R.~B.}} \AND
\bauthor{\bsnm{Kaminski},~\bfnm{D.}\binits{D.}}
(\byear{2001}).
\btitle{Asymptotics and {M}ellin--{B}arnes Integrals}.
\bseries{Encyclopedia of Mathematics and Its Applications}
\bvolume{85}.
\bpublisher{Cambridge Univ. Press}, \baddress{Cambridge}.
\bid{mr={1854469}}
\end{bbook}
%
\endbibitem

\bibitem{Patie-OU-06}
%
\begin{barticle}[mr]
\bauthor{\bsnm{Patie},~\bfnm{Pierre}\binits{P.}}
(\byear{2007}).
\btitle{Two-sided exit problem for a spectrally negative {$\alpha$}-stable
{O}rnstein--{U}hlenbeck process and the {W}right's generalized hypergeometric
functions}.
\bjournal{Electron. Comm. Probab.}
\bvolume{12}
\bpages{146--160 (electronic)}.
\bid{mr={2318162}}
\end{barticle}
%
\endbibitem

\bibitem{Patie-08a}
%
\begin{barticle}[mr]
\bauthor{\bsnm{Patie},~\bfnm{P.}\binits{P.}}
(\byear{2008}).
\btitle{{$q$}-invariant functions for some generalizations of the
{O}rnstein--{U}hlenbeck semigroup}.
\bjournal{ALEA Lat. Am. J. Probab. Math. Stat.}
\bvolume{4}
\bpages{31--43}.
\bid{mr={2383732}}
\end{barticle}
%
\endbibitem

\bibitem{Patie-CBI-09}
%
\begin{barticle}[mr]
\bauthor{\bsnm{Patie},~\bfnm{Pierre}\binits{P.}}
(\byear{2009}).
\btitle{Exponential functional of a new family of {L}\'evy processes and
self-similar continuous state branching processes with immigration}.
\bjournal{Bull. Sci. Math.}
\bvolume{133}
\bpages{355--382}.
\bid{doi={10.1016/j.bulsci.2008.10.001}, mr={2532690}}
\end{barticle}
%
\endbibitem

\bibitem{Patie-08-Asym}
%
\begin{barticle}[mr]
\bauthor{\bsnm{Patie},~\bfnm{P.}\binits{P.}}
(\byear{2009}).
\btitle{A few remarks on the supremum of stable processes}.
\bjournal{Statist. Probab. Lett.}
\bvolume{79}
\bpages{1125--1128}.
\bid{doi={10.1016/j.spl.2009.01.001}, mr={2510779}}
\end{barticle}
%
\endbibitem

\bibitem{Patie-09-cras}
%
\begin{barticle}[mr]
\bauthor{\bsnm{Patie},~\bfnm{Pierre}\binits{P.}}
(\byear{2009}).
\btitle{Law of the exponential functional of one-sided {L}\'evy
processes and
{A}sian options}.
\bjournal{C. R. Math. Acad. Sci. Paris}
\bvolume{347}
\bpages{407--411}.
\bid{doi={10.1016/j.crma.2009.02.013}, mr={2537239}}
\end{barticle}
%
\endbibitem

\bibitem{Patie-Asian-09}
%
\begin{bmisc}[auto:STB|2010-11-18|09:18:59]
\bauthor{\bsnm{Patie},~\bfnm{P.}\binits{P.}}
(\byear{2011}).
\bhowpublished{A Geman--Yor formula for one-sided {L}\'{e}vy processes.
Univ. Libre de Bruxelles. Unpublished manuscript}.
\end{bmisc}
%
\endbibitem

\bibitem{Patie-06c}
%
\begin{barticle}[mr]
\bauthor{\bsnm{Pierre},~\bfnm{Patie}\binits{P.}}
(\byear{2009}).
\btitle{Infinite divisibility of solutions to some self-similar
integro-differential equations and exponential functionals of {L}\'evy
processes}.
\bjournal{Ann. Inst. H. Poincar\'e Probab. Statist.}
\bvolume{45}
\bpages{667--684}.
\bid{doi={10.1214/08-AIHP182}, mr={2548498}}
\end{barticle}
%
\endbibitem

\bibitem{Rivero-05}
%
\begin{barticle}[mr]
\bauthor{\bsnm{Rivero},~\bfnm{V{\'{\i}}ctor}\binits{V.}}
(\byear{2005}).
\btitle{Recurrent extensions of self-similar {M}arkov processes and
{C}ram\'er's condition}.
\bjournal{Bernoulli}
\bvolume{11}
\bpages{471--509}.
\bid{doi={10.3150/bj/1120591185}, mr={2146891}}
\end{barticle}
%
\endbibitem

\bibitem{Titchmarsh-39}
%
\begin{bbook}[mr]
\bauthor{\bsnm{Titchmarsh},~\bfnm{E.~C.}\binits{E.~C.}}
(\byear{1939}).
\btitle{The Theory of Functions}, \bedition{2nd} ed.
\bpublisher{Oxford Univ. Press}, \baddress{London}.
\end{bbook}
%
\endbibitem

\bibitem{Widder-41}
%
\begin{bbook}[mr]
\bauthor{\bsnm{Widder},~\bfnm{David~Vernon}\binits{D.~V.}}
(\byear{1941}).
\btitle{The {L}aplace {T}ransform}.
\bseries{Princeton Mathematical Series}
\bvolume{6}.
\bpublisher{Princeton Univ. Press}, \baddress{Princeton, NJ}.
\bid{mr={0005923}}
\end{bbook}
%
\endbibitem

\bibitem{Yor-01}
%
\begin{bbook}[mr]
\bauthor{\bsnm{Yor},~\bfnm{Marc}\binits{M.}}
(\byear{2001}).
\btitle{Exponential Functionals of {B}rownian Motion and Related Processes}.
\bpublisher{Springer}, \baddress{Berlin}.
\bid{mr={1854494}}
\end{bbook}
%
\endbibitem

\end{thebibliography}
\end{document}